\documentclass[a4paper,10 pt,journal]{IEEEtran}  
\usepackage{mathrsfs, amsmath, amssymb, bm, bbm,bbold}
\usepackage{mathtools}
\usepackage{graphicx}
\usepackage[noadjust]{cite}
\usepackage{color}
 
\usepackage{gensymb}
\newtheorem{theorem}{Theorem}
\newtheorem{lemma}{Lemma}

\title{\LARGE \bf
A Sufficient Condition for Small-Signal Stability and Construction of Robust Stability Region}

\author{Parikshit Pareek$^{1}$, Konstantin Turitsyn$^{2}$, Krishnamurthy Dvijotham$^{3}$, and Hung D. Nguyen$^{4}$ 
}
\begin{document}

\maketitle
\thispagestyle{empty}
\pagestyle{empty}

\begin{abstract}
\color{black}
The small-signal stability is an integral part of the power system security analysis. The introduction of renewable source related uncertainties is making the stability assessment difficult as the equilibrium point is varying rapidly. This paper focuses on the Differential Algebraic Equation (DAE) formulation of power systems and bridges the gap between the conventional reduced system and the original one using logarithmic norm. We propose a sufficient condition for stability using Bilinear Matrix Inequality and its inner approximation as Linear Matrix Inequality. Another contribution is the construction of robust stability regions in state-space in contrast to most existing approaches trying same in the parameter space. Performance evaluation of the sufficiency condition and its inner approximation has been given with the necessary and sufficient condition for a small-scale test case. The paper provides a necessary base to develop tractable construction techniques for the robust stability region of power systems. 
\color{black}
\end{abstract}

\section{INTRODUCTION}

The small-signal stability analysis 
is used to measure and identify electrical power system's capacity to remain in  synchronism under the influence of small disturbance. The state matrix provides information about the stability as well as about oscillations and different modes of oscillation in power systems \cite{kundur1994power}. 

The behavior of eigenvalue with different state and parameter of the power system has been the focus of the small-signal stability studies. The stability is assessed with variations in operating point, leading to a new equilibrium point, using the necessary and sufficient condition of stability on the eigenvalue of state matrix. Similarly, the parameter variation leading to change in the dynamics of the system as also been studied. The graphical method of root locus analysis, requiring repeated calculations, has also been used to track the eigenvalue movement with a change in state or parameter \cite{ogata2002modern,chow2004power}.   

The movement of eigenvalues with the change in the base point has been explored in many works related to eigenvalue sensitivity \cite{pagola1989sensitivities,luo2011sensitivity}. The continuation of invariant subspace (CIS) method has been proposed to calculate the successive eigenvalue sensitivities \cite{yang2007critical}. The sensitivity is a function of state and parameters of the system and keeps changing with variations. The validity range of sensitivity cannot be specified explicitly. Therefore, the eigenvalue information is local and hence cannot be used to identify the stable region. 

Recently, there are some attempts to identify the small signal boundaries. In \cite{li2013novel}, authors have presented a method which can be used for the stability boundary calculation in the event of variation in the power system stabilizer (PSS) time constants. The process is based on the concepts of eigenvalue perturbation and sensitivity. On similar lines, a strategy for computing small-signal security margins has been presented recently, computing stability margin for load increase, in \cite{li2018investigation}. A method based on Gersgorin disks has also been used to identify the region of eigenvalue perturbation for disturbances in micro-grid system \cite{li2016gervsgorin}. An extension using reachable set computation (FAR) has also been presented in  \cite{li2018formal}. Both these methods work in complex eigenvalue space and need to calculate the disturbance matrix each time. 

A method to construct robust stability region in parameter space for PSS has been presented recently in \cite{pulcherio2018robust}. The work uses the bialternate sum matrix approach and shown to be efficient over the branch and bound approach used in \cite{pai1997applications} for similar purpose. This approach works only if the reduced system matrix is affine and hence works for parameters variations of the system but not with the state variations.  

This work is a follow-up to the contraction analysis presented in \cite{nguyen2016contraction}. In this paper, a relation between the reduced form and the associated DAE system has been established using logarithmic norm. Further, we use this relation to develop a sufficient condition for stability of the DAE system. The robust stability region on which, if there exists any equilibrium point then, the system will be small-signal stable at that equilibrium, has been constructed. The sufficient condition for stability, as a Bilinear Matrix Inequality (BMI) is derived using the logarithmic norm. Then, a convex inner approximation of the BMI has been presented as Linear Matrix Inequality (LMI). A 2-bus test system example is used to present the simulation results of the sufficient stability conditions. The local validity of the eigenvalue and its sensitivity concepts has also been shown. The results of the robust stable area with sufficient BMI and its convex approximation as LMI have also been presented in state space. 

\color{black}
Main contributions of the paper are as follows.
\begin{enumerate}
   \item Development of a relationship between the reduced system and DAE system using logarithmic norm. 
    \item A sufficient condition for small-signal stability has been developed for the DAE form.
    \item Performance evaluation of sufficient condition and its inner approximation has been given with the necessary and sufficient condition for a small-scale test case.  
\end{enumerate}
\color{black}

\section{Main Result}\label{sec:mainresults}
In this paper, we use $S^+_n$ to denote the set of symmetric positive definite matrices of dimension $n$. $\lambda(A)$ denotes eigenvalues of matrix $A \in \mathbb{R}^{n \times n}$.

For a power system, there are two sets of equations describing its behavior. The differential equation set, $\mathbf{\dot x}=f(\mathbf{x},\mathbf{y})$, is used to describe the dynamics of power system. The prime mover-generator system dynamics is described using the swing equation, which is a relation between acceleration in load angle and accelerating power. Other equations in this set involve the relations of internal voltage change with generator parameters and axis component of voltage and current. The number of differential equations depends upon the modeling of the system and components modeled along with the generator.  The time constants involved in these are much higher than the time constant of electrical dynamic equations of transmission components like transmission lines and transformers. Therefore, the fast dynamics are ignored and the transmission system coupled with generator stator is modeled as a set of algebraic equations, $ {0}=g(\mathbf{x},\mathbf{y})$. Both of these sets of equations together are also called the ODEs with constraints. In both sets, the vector $\mathbf{x}\in \mathbb{R}^n$ contains dynamic state variables while $\mathbf{y}\in \mathbb{R}^m$ represents the algebraic variables. Therefore, the system dynamics can be expressed as semi-explicit DAE as: 
\begin{equation}
\begin{aligned}\label{eq:DAE}
    \mathbf{\dot x} &= f(\mathbf{x},\mathbf{y}) \\
             {0}    &= g(\mathbf{x},\mathbf{y})
\end{aligned} 
\end{equation}

We further assume that, for the considered states $\mathbf{x}$, there exists corresponding variables $\mathbf{y}$ satisfying the algebraic constraints. For stability studies, general practice is the elimination of the algebraic variables via reduction techniques. This elimination will leave us with a modified set of ODEs $(\mathbf{\dot x}=f_{1}(\mathbf{x}))$. This method disallows to explore the structural properties of the DAE form. To this reason, a variety of literature focus on the DAE form, \cite{wu1994stability,lewis1985optimal} especially in the context of descriptor system $ E \,  \mathbf{\dot z} = h (\mathbf{z})$; where, $\mathbf{z^T}=[\mathbf{x^T}, \mathbf{y^T}]$, $h^T=[f^T,g^T]$ and $E$ is a diagonal matrix of $ \mathbb{R}^{(n+m)\times(n+m)}$ with $E_{ii}=1$ if $i\leq n$ and zero otherwise. 
In the power system, the structure-preserving DAE model is gaining importance in eigenvalue-based stability analysis \cite{li2017network},\cite{li2018investigation}  or for contraction region \cite{nguyen2016contraction}. To obtain the power system equivalent of the DAE system (\ref{eq:DAE}) we define a block Jacobian as:
\begin{equation}\label{eq:jacobian}
J(\mathbf{x},\mathbf{y})=
    \begin{bmatrix}
{\partial f}/{\partial \mathbf{x}} & {\partial f}/{\partial \mathbf{y}}\\ {\partial g}/{\partial \mathbf{x}}& {\partial g}/{\partial \mathbf{y}}  
    \end{bmatrix}=
    \begin{bmatrix}
A & B\\ C& D
    \end{bmatrix}.
\end{equation}

In simplifying notations, the linearized state space model with algebraic constraints will look like:
\begin{align}\label{eq:ssPS}
    \delta\mathbf{\dot x} &= A \delta\mathbf{x}+B \delta\mathbf{y},\\
    0 & =C \delta\mathbf{x} + D \delta\mathbf{y}.
\end{align}

To assess the stability of the above system, common practice is to eliminate the algebraic variables $\delta{y}$ and reduce the system to $\delta\mathbf{\dot x} = J_r \delta\mathbf{x}$ where $J_r = A -B D^{-1}C$ be the reduced Jacobian matrix. The reduced system is small-signal stable if and only if all eigenvalues of the reduced Jacobian lie in the left half plane or $\mathrm{Re}\left(\lambda(J_r)\right)<0$ \cite{chen1998linear}.


\begin{table}[ht]
  \centering
  \caption{Standard matrix measures}
    \begin{tabular}{|c|c|}
    \hline
    Vector norm, $ \|\cdot\|$ & Matrix measure, $\mu_p(M)$ \\ [4.5pt]
    \hline
   $\| x\|_1 = \sum_{i}|x_i|$ & $\mu_1(M)=\max_j(m_{jj}+\sum_{i\neq j}|m_{ij}|)$ \\ [4.5pt]
$\|x\|_2 = (\sum_{i}|x_i|^2)^{1/2}$ & $\mu_2(M)=\max_i(\lambda_i\{\frac{M+M^T}{2}\})$ \\[4.5pt] 
  $\|x\|_\infty = \max_{i}|x_i|$ & $\mu_\infty(M)=\max_i(m_{ii}+\sum_{j\neq i}|m_{ij}|)$\\[4.5pt]
  \hline
    \end{tabular}%
  \label{tab:matmu}%
\end{table}
Before introducing the main result, we define the logarithmic norm or the matrix measure. The matrix measure $\mu_p(M)$ of a matrix $M$ is defined as \cite{vidyasagar2002nonlinear}:
\begin{equation}
    \mu_p(M) := \lim\limits_{\epsilon\rightarrow0^+}\frac{1}{\epsilon}(\|I+\epsilon M\|_p - 1)
\end{equation}

Unlike norms, $\mu_p(M)$ can be negative as well. Thus, it can be an important yardstick to construct and measure the stable region providing upper bound on the real part of eigenvalues of a matrix $M$ and relating stability of ODE and DAE. For $p=1,\,2,\,\infty$, the standard matrix measures as well as vector norms are listed in Table \ref{tab:matmu}.

We define \textit{generalized reduced Jacobian} $F_r = P J_r$ with $P \in S^+_n$, and \textit{generalized unreduced Jacobian} $F = Z^T J$ which can be written as 
\begin{align} 
   F &= \begin{bmatrix} PA + R^T C  &  PB + R^T D \\
   Q^T C & Q^T D\end{bmatrix}.\label{eq:F}
   \end{align}
Here $Z= \begin{bmatrix} P & \mathbf{0} \\ R & Q \end{bmatrix}$ is an auxiliary matrix consisting of $\mathbf{0}\in \mathbb{R}^{n\times m}$ a null matrix, $R\in \mathbb{R}^{m\times n}$, and $Q\in \mathbb{R}^{m\times m}$. The main result is presented below.
\begin{lemma} \label{lem:main}
For the DAE system \eqref{eq:DAE} with a negative logarithmic norm $\mu_p(F) < 0$ characterized by the matrices $P,Q,R$, the following relation holds:
\begin{equation} \label{eq:FrF}
	\mu_p(F_r) \leq \mu_p(F).
\end{equation}
\end{lemma}

The proof follows closely the reasoning to prove Theorem 1 for contraction analysis in \cite{nguyen2016contraction}; however, the generalized reduced and unreduced Jacobian are slightly different from what we define in this paper. Note also that here we focus on the original variables $\delta \mathbf{z} = (\delta \mathbf{x}^T, \delta \mathbf{y}^T)^T$ instead of their transformed counterparts $\delta \mathbf{w} = (\delta \mathbf{v}^T, \delta \mathbf{u}^T)^T$ in contraction analysis. Using Lemma \ref{lem:main}, we introduce a sufficient condition of stability in terms of logarithmic norm. 

\begin{theorem} \label{theo:main}
The DAE system is small-signal stable if there exists a generalized unreduced Jacobian has a negative matrix measure, i.e., $\mu_p(F)<0$. 
\end{theorem}
Let consider a Lyapunov candidate function $V = P \|\delta \mathbf{x}\|_p$ with $P \in S^+_n$. Following the same reasoning for matrix measure results introduced in \cite{vidyasagar2002nonlinear, aminzarey2014contraction}, one takes upper Dini derivative of $V$ to yield:
\begin{equation} \label{eq:DV}
D^+\|\delta \mathbf{x}(t)\|_p 
\leq \mu_p (F_r)\, \|\delta \mathbf{x}(t)\|_p.
\end{equation}

Combining \eqref{eq:DV} and Lemma \ref{lem:main}, a negative $\mu_p(F)$ verifies that $V$ is a valid Lyapunov function, thus the system is small-signal stable at the considered equilibirum point.

While most of the stability assessment approaches presented in the past rely on the reduced form, Lemma \ref{lem:main} connects it to the DAE system. However, Theorem \ref{theo:main} can be proved without Lemma \ref{lem:main}, by defining another Lyapunov function $V_1 = \|Z^T \delta \mathbf{z}\|_p$ and showing that $D^+V_1 
\leq \mu_p (F)\, V_1$.

\color{black}
\subsection*{Special case for $p=2$:}
As most of the existing work focuses on eigenvalue analysis and quadratic Lyapunov functions which relate to norm $2$, in the following, we reproduce the main result in Lemma \ref{lem:main} for $p=2$, by using \textit{Interlacing theorem}. 

From the definitions of standard logarithmic norms given in Table \ref{tab:matmu}, $\mu_2(F)=\lambda_{max}((F+F^T)/2)$. Therefore, the condition $\mu_2(F)< 0$ implies that $\lambda_{max}(F+F^T)< 0$.
By defining $R=-(PBD^{-1})^T$ in definitions of $F$ and $F_r$ from Theorem \ref{lem:main}, $F+F^T$ will become:
\begin{equation}\label{eq:L1}
 F+F^T=   \begin{bmatrix} F_r+Fr^T & C^TQ\\ Q^TC & Q^TD+D^TQ
    \end{bmatrix}.
\end{equation}

Now, as $F+F^T$ is a real symmetric matrix, with ordered eigenvalue sequence of $F+F^T$ as $\lambda_1\leq \lambda_2 \leq \dots \leq \lambda_{n+m}$, by \textit{Interlacing Theorem} \cite{smith1992some} for Hermittian matrices (for $i=1\dots n$):
\begin{equation*}
    \lambda_i(F_r+F_r^T)  \leq \lambda_{i+m}(F+F^T)
\end{equation*}
which leads to 
\begin{equation} \label{eq:interlacing}
    \lambda_{max}(F_r+F_r^T)  \leq \lambda_{max}(F+F^T).
\end{equation}

Combining (\ref{eq:interlacing}) with the definition of log norm ($\mu_2$), it is proved that $\mu_2(F_r) \leq \mu_2(F)$, which is Lemma \ref{lem:main} for $p=2$. In the following sections, we use $\mu_2(F)$ to construct the stability region. 
\color{black}
\section{Construction of Robust Stable Region}\label{sec:region}
In this section, we describe the method to construct robust stability region within which the system is small-signal stable for every equilibrium point. The resulting region can be further incorporated in security-constrained optimization to provide convex search space \cite{nguyencons18, nguyenframe18, li2015tight, molzahn2018towards, golestaneh2018ellipsoidal,Jwang18}. We work on the Jacobian $J$, which can be expressed as an affine function of the system variables. This is obtained by expressing the DAEs in quadratic form of variables $\mathbf{z}$ \cite{nguyen2016contraction}. Thus, the system can be expressed as $J(\mathbf{z})=J_0+\sum_kz_kJ_k$ for $k=1 \dots (n+m)$.

\subsection{Stable Region by Bilinear Matrix Inequalities}

From Theorem \ref{theo:main}, the sufficient condition for stability of the system (\ref{eq:ssPS}) has been transformed from $\mathrm{Re}\left(\lambda(J_r)\right)<0$ to $\mu_p(F)< 0$. In view of the special case of $\mu_2$, the condition for stability can be stated as $\lambda_{max}(F+F^T)<0$. Therefore, the sufficient condition for stability can be expressed as a BMI in $Z$ and $J(\mathbf{z})$, with $\zeta$ as maximum eigenvalue and $I$ representing an identity matrix of size $(n+m)\times (n+m)$, as:
\begin{align}\label{eq:BMI}
     F(\mathbf{z},Z)+ F^T(\mathbf{z},Z)- \zeta I& \preceq 0.
\end{align}

Adding the conditions $\zeta<0$ and $P \succ 0$, the complete problem will involve a BMI with a set of constraints as LMIs. The BMI (\ref{eq:BMI}) alone turns out to be $\mathcal{NP}\textit{-hard}$ for solvability \cite{toker1995np} and define an $(n+m)- dimensional$ non-convex feasibility space. Therefore, there is no off-the-shelf algorithm to solve and construct the stability region using this BMI (\ref{eq:BMI}). Hence, it imperative to construct a convex inner approximation of the region defined by the BMI.  

\subsection{Inner Approximation of BMI as Robust LMI}
To convert the bilinear relation into the linear one, we fix $Z=Z^\star$ in the BMI (\ref{eq:BMI}). The resultant defines a convex region of small-signal stability expressed as an LMI in variable $\mathbf{z}$. As the solution space of an LMI is convex, fixing of one variable will provide the convex inner approximation of the BMI. The construction of $Z^\star$ is done based upon the sufficient condition proved in Theorem 1 using the Jacobian $J(\mathbf{z^*})$ at the equilibrium point. The set of LMIs with variable matrices $P$ and $Z$ is:
\begin{equation}\label{eq:LMI}
    \begin{aligned}
   &&  P \succ 0,\\
   && J(\mathbf{z^*})^TZ+Z^TJ(\mathbf{z^*})-\zeta I \preceq 0.\\
        \end{aligned}
\end{equation}

A feasible solution of this LMI set, with constraint $\zeta < 0$, will give the matrix $Z^\star$. Now, the inner approximated convex stability region of the non-convex region defined by the BMI (\ref{eq:BMI}) can be obtained by solving a robust SDP problem for affine perturbation for $\mathbf{z}\in \mathcal{U}$:
\begin{equation}\label{eq:RSDP}
\begin{aligned}
    & {\text{minimize}} & & \zeta \\
& \text{subject to}
& & J(\mathbf{z})^T Z^\star+(Z^\star)^T J(\mathbf{z})-\zeta I \preceq 0
\end{aligned}
\end{equation}

If the solution $\zeta^\star$ of this robust SDP is negative, then the region defined by the set $\mathcal{U}$ is the robust stability region. As most of the operational constraints are box-typed, we restrict our self to a box-typed description of stability region in $z_k$. The box $\mathcal{U}$ is constructed by the set of $\mathcal{U}_k$ and $\mathcal{U}_k=[\underline{z_k}, \overline{z_k}]$. By taking leverage of sufficiency condition, robust SDP problem (\ref{eq:RSDP}) can be converted into a feasibility problem by adding a constraint as $\zeta < 0$. 

\subsection{Eigenvalue Sensitivity}
The information provided by eigenvalue sensitivity has been used widely to examine various aspects of stability. The analytic expression of sensitivity eases out the computation and provides a better understanding about the variation. The matrix pencil ($J(\mathbf{z}),E$) has two sets of eigenvalues at any particular $\mathbf{z}$. One with $n$ finite eigenvalues and another one with $m$ infinite eigenvalues. It can be shown that all $n$ finite eigenvalues of ($J(\mathbf{z}),E$) are same as $n$ eigenvalues of $J_r$ \cite{li2013novel}. As the Jacobian (\ref{eq:jacobian}) is an affine function of state variables, the analytic expression for eigenvalue sensitivity of ($J(\mathbf{z}),E$) with respect to $\mathbf{z}$ can be developed, which gives sensitivity equivalent to that for eigenvalues of $J_r$.

If $\lambda$ is an eigenvalue of the matrix pencil ($J(\mathbf{z}),E$) with left and right eigen vectors as $\upsilon$ and $\nu$, then with similar approach presented in \cite{li2013novel}, we obtain:
\begin{equation}\label{eq:dlamda}
\frac{\partial \lambda_i}{\partial\mathbf{z}_k} =\frac{\upsilon^T_iJ_k\nu_i}{\upsilon^T_iE\nu_i}.
\end{equation}

This sensitivity expression is used to show the local validity of sensitivity based approaches. This will be further utilized to construct stability regions in future works. 


\color{black}

\section{Simulation}
\subsection{$2$-bus Test System}\label{sec:2bus}
\begin{figure}[t]
    \centering
    \includegraphics[width=0.7\columnwidth]{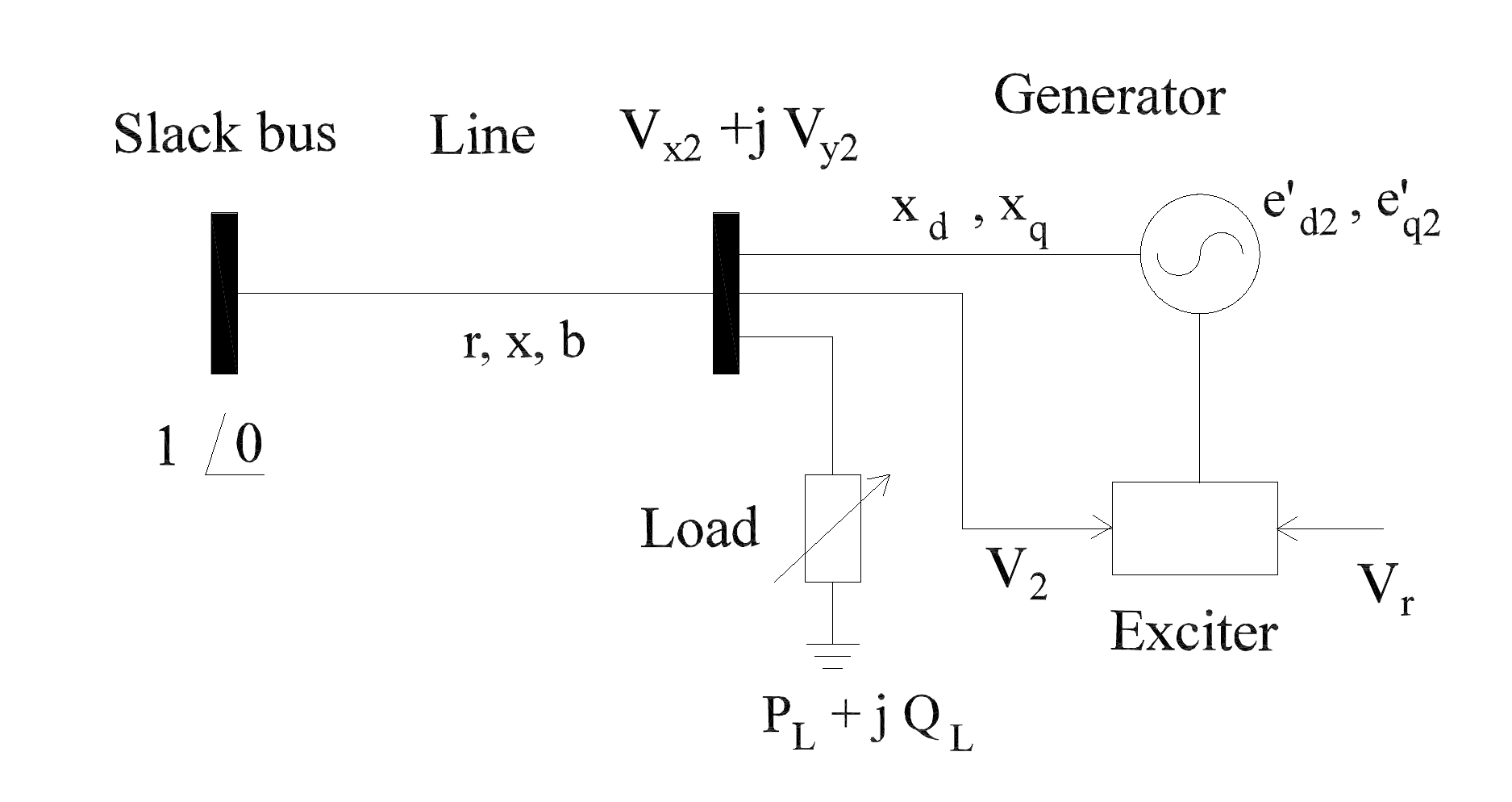}
	\caption{A $2$-bus test system \cite{nguyen2016contraction}}
    \label{fig:rudi}
\end{figure}

The construction approaches of robust stability region as BMI and LMI along with analytic expression of eigenvalue sensitivity are discussed above. In this section, we illustrate the procedure by constructing the stability region for a two-bus system as shown in Figure \ref{fig:rudi} \cite{nguyen2016contraction}.

A PSAT dynamic model \cite{PSAT} input based dynamic simulation and PowerWolf, an analysis package developed in Mathematica $11.1.0.0$ by Turitsyn's group at MIT are used for system modeling. We also use YALMIP with Sedumi \cite{lofberg2012} in MATLAB for solving robust LMIs.

\color{black}
\subsection{Results and Analysis} 
\begin{figure}[b]
    \centering
    \includegraphics[width=0.9\columnwidth]{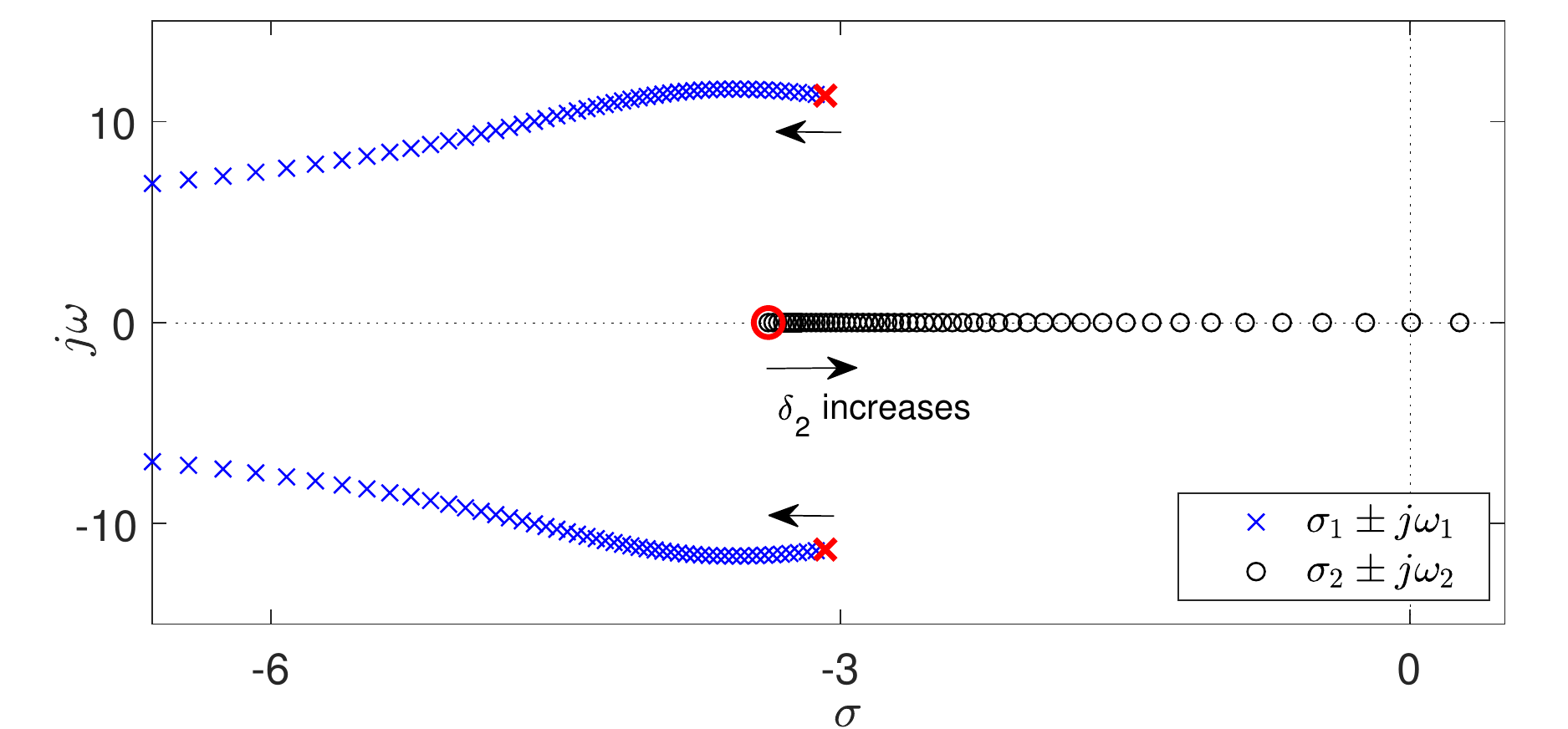}
	\caption{Root locus plot with increase in $\delta_2$}
	\label{fig:root_locus}
\end{figure}
\begin{figure}[t]
    \centering
    \includegraphics[width=0.95\columnwidth]{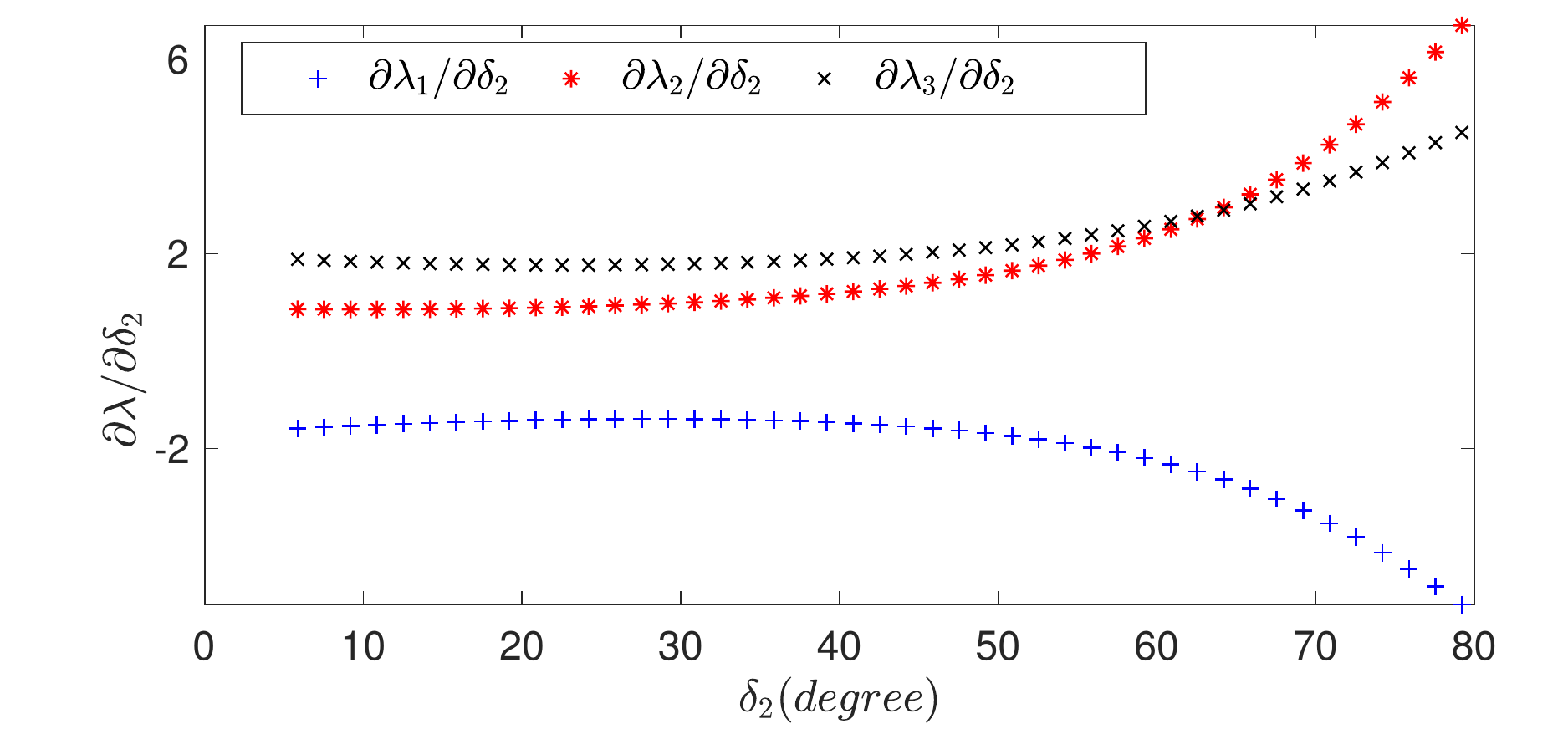}
	\caption{Variation of the eigenvalue sensitivity $\partial\lambda/\partial\delta_2$ with $\delta_2$}
	\label{fig:sensitivity}
\end{figure}
\begin{figure}[t]
    \centering
    \includegraphics[width=0.95\columnwidth]{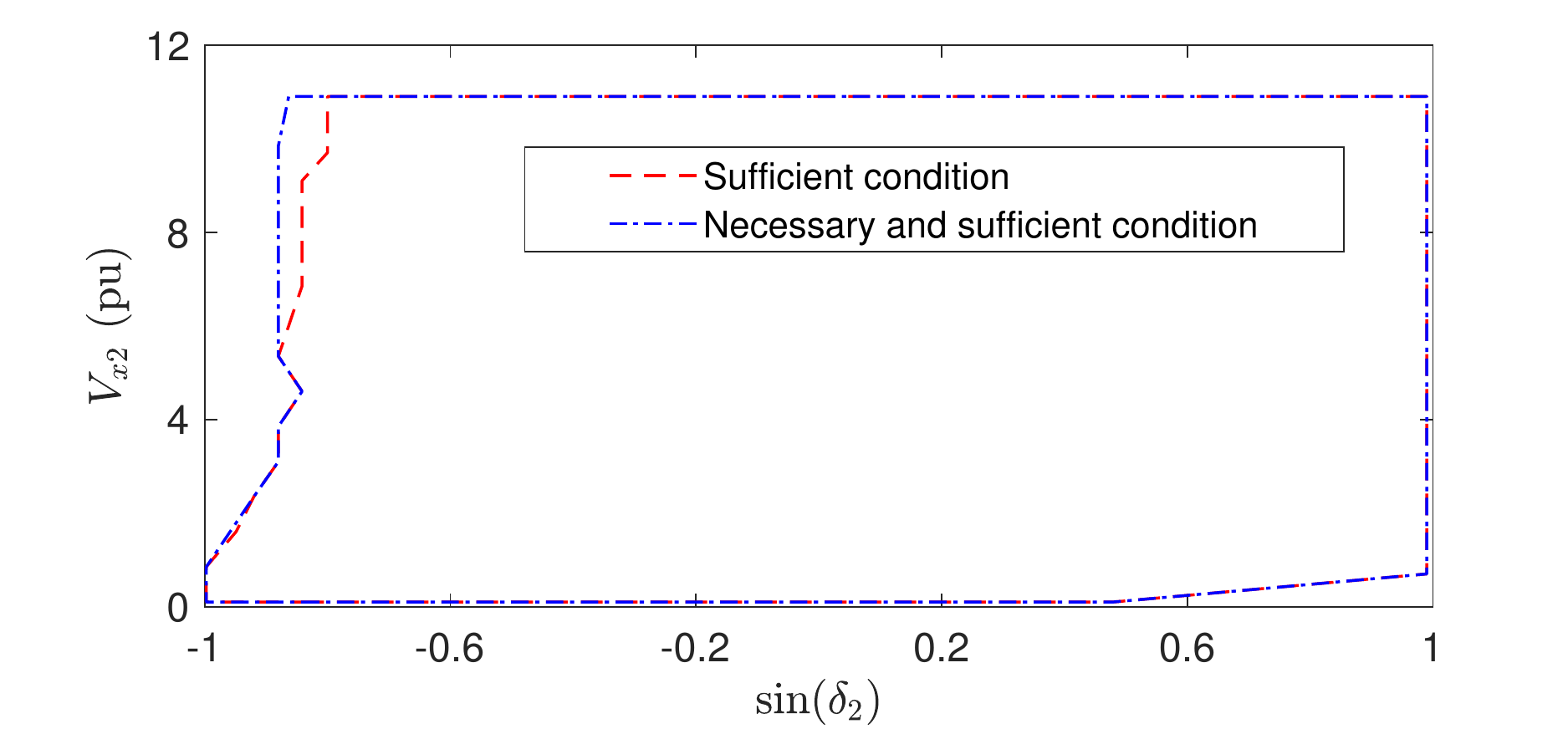}
	\caption{Performance of sufficient condition (BMI \ref{eq:BMI}) in 2-D state plane}
	\label{fig:BMI area}
\end{figure}
\begin{figure}[h]
    \centering
    \includegraphics[width=0.9\columnwidth]{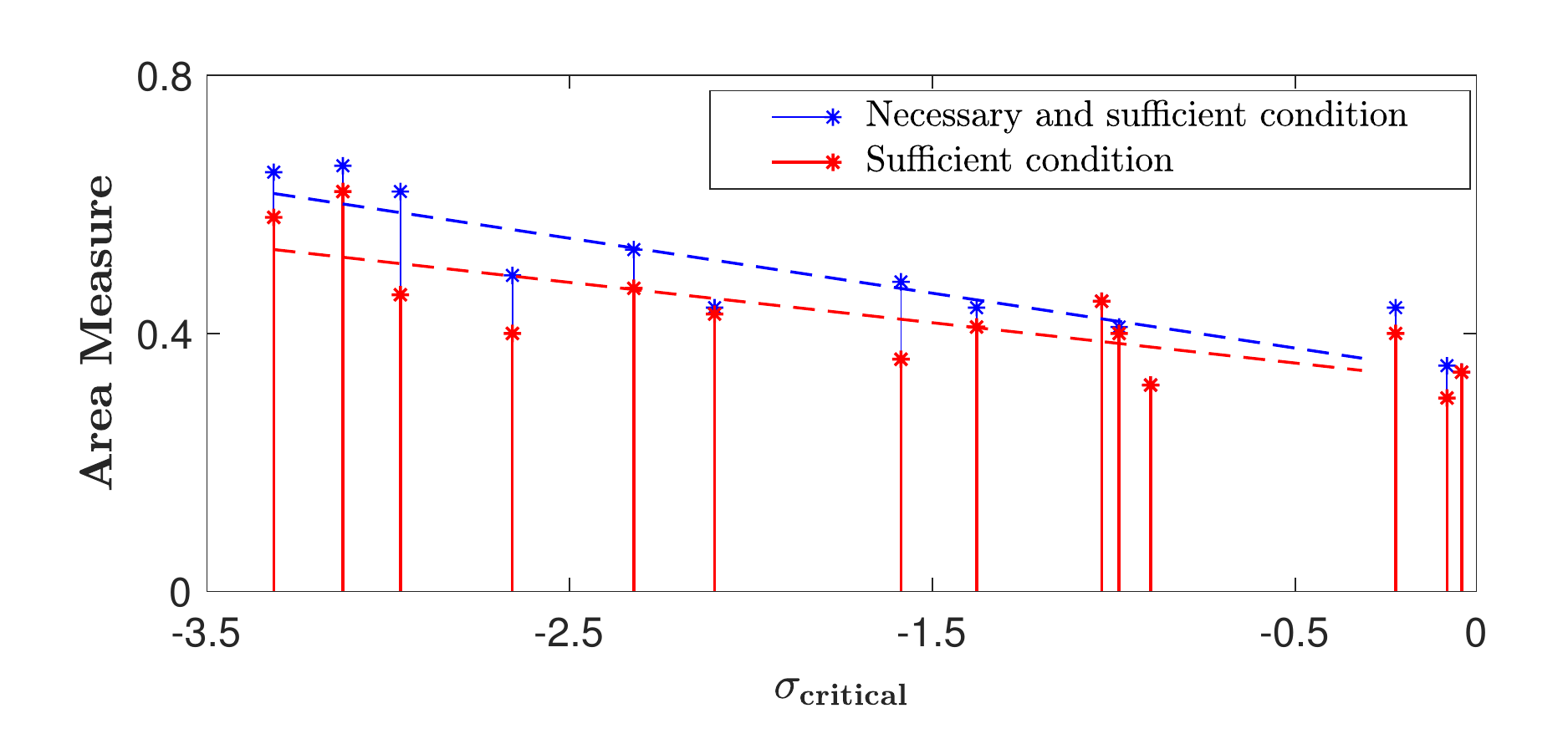}
	\caption{Stability area measure by different stability conditions}
	\label{fig:area}
\end{figure}
We begin with showing the limitations of eigenvalue approaches to underline the importance of developed sufficient condition. Figure \ref{fig:root_locus} is root-locus plot showing the movement of eigenvalues ($\lambda_i=\sigma_i+j\omega_i$) of matrix pencil ($J(\mathbf{z}),E$) with increase in $\delta_2$. The critical most eigenvalue at base point is a complex-conjugate pair with $\sigma=-3.07$ shown as big red cross. As the value of $\delta_2$ increases, this eigenvalue pair moves away from the $j\omega$-axis while the other eigenvalue, with real part only, moves towards the origin. The $\sigma_2$ crosses $\sigma_1$ at $\delta_2=10.06\degree$ and becomes the critical most eigenvalue. This shows that critical eigenvalue is state dependent. Hence, any region calculation based upon critical eigenvalue at the current state can lead to overestimation. Figure \ref{fig:root_locus} also shows that the spacing between two successive values of $\lambda$ is not constant with constant $\Delta\delta_2$. This is because of variable eigenvalue sensitivity, shown in Figure \ref{fig:sensitivity}. The sensitivity variation clearly proves that $\partial\mathbb{\lambda}_{i}/\partial\delta_2$ is not constant. Thus, the sensitivity matrix requires continuous updates and is not very suitable for the stability region identification in state space. 

In Figure \ref{fig:BMI area}, the stable area is plotted by the necessary and sufficient condition of stability, $\sigma_{critical}(J_r(\mathbf{z}))< 0$ and the sufficient condition (\ref{eq:BMI}), in $2-D$ state plane. The stable boundary is obtained via verifying stability of a very large number of points in the plane. We limit ourself to $2$-bus test system as construction of stable space is very expensive in this manner. The performance of the proposed sufficiency condition has been proved to be impressive by this plot.

As both stability conditions are non-convex in $n+m-dimension$ state plane, an area measure is defined to quantify the stable area as the ratio of stable points to the total number of points generated and tested for the stability in the $4-D$ space. The Figure \ref{fig:area} shows that the area by BMI is always less than equal to the area by the necessary and sufficient condition. Dotted lines indicate the single variable linear regression relation between area and $\sigma_{critical}$. Negative slope suggests that the stable area tend to decrease when $\sigma_{critical}$ moves towards the zero. Also, the decreasing distance between two regression lines suggest that sufficiency condition (BMI \ref{eq:BMI}) tend to improve as $\sigma_{critical}$ approaches to the origin.

\begin{figure}[t]
    \centering
    \includegraphics[width=0.9\columnwidth]{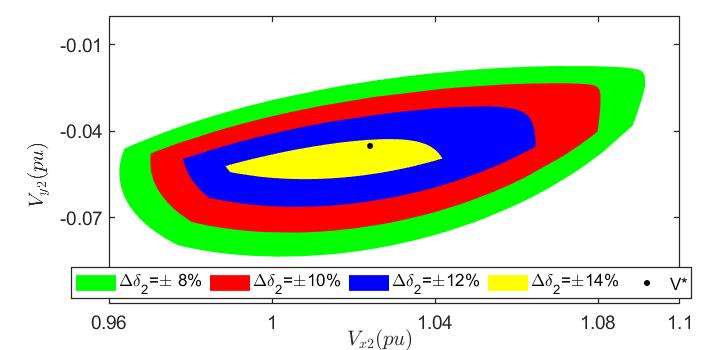}
	\caption{Decrease in stable space in $V_{x2}-V_{y2}$ plane with increase in $\Delta\delta_2$}
	\label{fig:Vx_Vy_space}
\end{figure}

\begin{figure}[t]
    \centering
    \includegraphics[width=0.9\columnwidth]{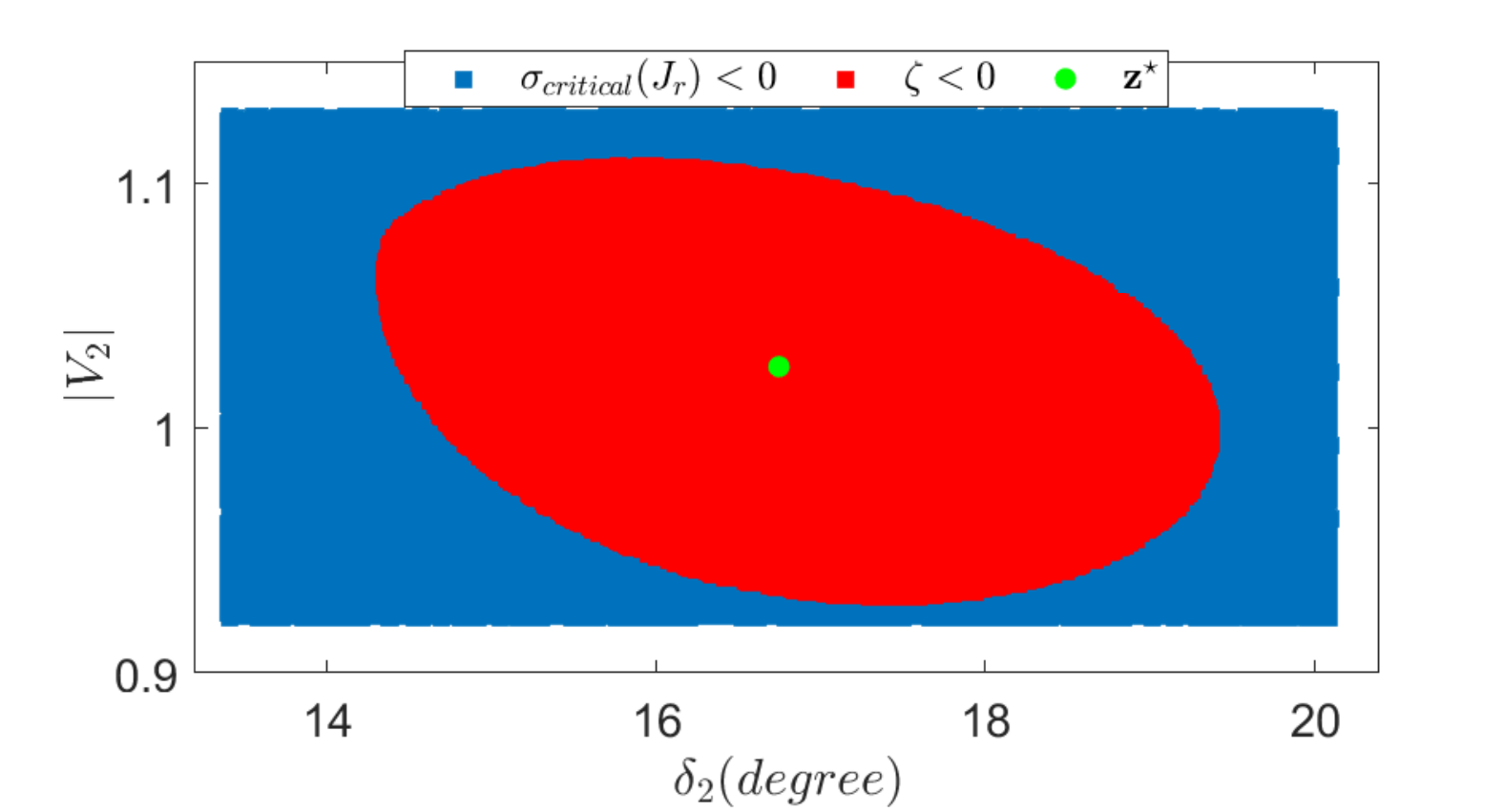}
	\caption{Stable space in $\delta_2-|V_2|$ plane using LMI (\ref{eq:RSDP})}
	\label{fig:Del-V space}
\end{figure}

By solving robust SDP (\ref{eq:RSDP}) for feasibility, Figure \ref{fig:Vx_Vy_space} is plotted for different $\mathcal{U}_k$ for variable $\sin{\delta_2}$ and $\cos{\delta_2}$. The uncertainty space is constructed using box-type uncertainty in the system variable $\delta_2$ and then it gets converted into elliptical one due to the algebraic relationship. Similar relationship can be observed between the variable space transformation from $V_{x2}-V_{y2}$ variable to $|{V_2}|$. The result shown in Figure \ref{fig:Vx_Vy_space} highlights the dependency of stability space in different planes with each other. For example the red space indicates all the values of $V_{x2}$ and $V_{y2}$ for which system is stable if $\Delta\delta_2$ attains any value between $\pm 10\%$. The robust stability region, around the base point, is a rectangle made up with  $\mathcal{U}_{V_{x2}}=\mathcal{U}_{V_{y2}}$ defined by $\pm 5\%$ variation in the corresponding states and $\mathcal{U}_{\sin{\delta_{2}}}=\mathcal{U}_{\cos{\delta_{2}}}$ defined by $\pm 8\%$. For variations of $\pm 10\%$ in $|V_2|$ and $\pm 20\%$ in $\delta_2$, the inner approximated robust stable space is shown in Figure \ref{fig:Del-V space}. The blue points cover entire rectangular area while the red ellipse is superimposed over it. This is the feasible solution space of sufficient condition expressed as robust SDP (\ref{eq:RSDP}) and thus convex in nature. This space provides stability certificate for all equilibrium points, within it, without actually constructing state matrix and evaluating the eigenvalues at those points. 
\section{CONCLUSIONS}\label{sec:conclusion}  
This paper deals with the construction and analysis of the robust small-signal stability region. We provide a sufficient condition for stability as BMI along with LMI approximation of the same. The paper offers a path for developing efficient methods for robust stability region identification for large-scale power systems using LMIs. We further, intend to look deeper into the shape of the region using the sensitivity information. The generalization of the stable region geometry using the $n-ellipse$ with multiple foci will also be explored. 
\color{black}
\bibliographystyle{IEEEtran}
\bibliography{main.bbl}

\begin{thebibliography}{10}
\providecommand{\url}[1]{#1}
\csname url@samestyle\endcsname
\providecommand{\newblock}{\relax}
\providecommand{\bibinfo}[2]{#2}
\providecommand{\BIBentrySTDinterwordspacing}{\spaceskip=0pt\relax}
\providecommand{\BIBentryALTinterwordstretchfactor}{4}
\providecommand{\BIBentryALTinterwordspacing}{\spaceskip=\fontdimen2\font plus
\BIBentryALTinterwordstretchfactor\fontdimen3\font minus
  \fontdimen4\font\relax}
\providecommand{\BIBforeignlanguage}[2]{{%
\expandafter\ifx\csname l@#1\endcsname\relax
\typeout{** WARNING: IEEEtran.bst: No hyphenation pattern has been}%
\typeout{** loaded for the language `#1'. Using the pattern for}%
\typeout{** the default language instead.}%
\else
\language=\csname l@#1\endcsname
\fi
#2}}
\providecommand{\BIBdecl}{\relax}
\BIBdecl

\bibitem{kundur1994power}
P.~Kundur, N.~J. Balu, and M.~G. Lauby, \emph{Power system stability and
  control}.\hskip 1em plus 0.5em minus 0.4em\relax McGraw-hill New York, 1994,
  vol.~7.

\bibitem{ogata2002modern}
K.~Ogata and Y.~Yang, \emph{Modern control engineering}.\hskip 1em plus 0.5em
  minus 0.4em\relax Prentice hall India, 2002, vol.~4.

\bibitem{chow2004power}
J.~H. Chow, G.~E. Boukarim, and A.~Murdoch, ``Power system stabilizers as
  undergraduate control design projects,'' \emph{IEEE Transactions on power
  systems}, vol.~19, no.~1, pp. 144--151, 2004.

\bibitem{pagola1989sensitivities}
F.~L. Pagola, I.~J. Perez-Arriaga, and G.~C. Verghese, ``On sensitivities,
  residues and participations: applications to oscillatory stability analysis
  and control,'' \emph{IEEE Transactions on Power Systems}, vol.~4, no.~1, pp.
  278--285, 1989.

\bibitem{luo2011sensitivity}
C.~Luo and V.~Ajjarapu, ``Sensitivity-based efficient identification of
  oscillatory stability margin and damping margin using continuation of
  invariant subspaces,'' \emph{IEEE Transactions on Power Systems}, vol.~26,
  no.~3, pp. 1484--1492, 2011.

\bibitem{yang2007critical}
D.~Yang and V.~Ajjarapu, ``Critical eigenvalues tracing for power system
  analysis via continuation of invariant subspaces and projected arnoldi
  method,'' \emph{IEEE Transactions on power Systems}, vol.~22, no.~1, pp.
  324--332, 2007.

\bibitem{li2013novel}
C.~Li and Z.~Du, ``A novel method for computing small-signal stability
  boundaries of large-scale power systems,'' \emph{IEEE Trans. Power Syst.},
  vol.~28, no.~2, pp. 877--883, 2013.

\bibitem{li2018investigation}
C.~Li, H.-D. Chiang, and Z.~Du, ``Investigation of an effective strategy for
  computing small-signal security margins,'' \emph{IEEE Transactions on Power
  Systems}, 2018.

\bibitem{li2016gervsgorin}
Y.~Li, P.~Zhang, L.~Ren, and T.~Orekan, ``A ger{\v{s}}gorin theory for robust
  microgrid stability analysis,'' in \emph{PESGM, 2016}.\hskip 1em plus 0.5em
  minus 0.4em\relax IEEE, 2016, pp. 1--5.

\bibitem{li2018formal}
Y.~Li, P.~Zhang, and P.~B. Luh, ``Formal analysis of networked microgrids
  dynamics,'' \emph{IEEE Trans. on PS}, vol.~33, no.~3, pp. 3418--3427, 2018.

\bibitem{pulcherio2018robust}
M.~C. Pulcherio, M.~S. Illindala, and R.~K. Yedavalli, ``Robust stability
  region of a microgrid under parametric uncertainty using bialternate sum
  matrix approach,'' \emph{IEEE Transactions on Power Systems}, 2018.

\bibitem{pai1997applications}
M.~Pai, C.~Vournas, A.~Michel, and H.~Ye, ``Applications of interval matrices
  in power system stabilizer design,'' \emph{International Journal of
  Electrical Power \& Energy Systems}, vol.~19, no.~3, pp. 179--184, 1997.

\bibitem{nguyen2016contraction}
H.~D. Nguyen, T.~L. Vu, J.-J. Slotine, and K.~Turitsyn, ``Contraction analysis
  of nonlinear dae systems,'' \emph{arXiv preprint arXiv:1702.07421}.

\bibitem{wu1994stability}
H.~Wu and K.~Mizukami, ``Stability and robust stabilization of nonlinear
  descriptor systems with uncertainties,'' in \emph{IEEE CDC, 1994.}, vol.~3,
  1994, pp. 2772--2777.

\bibitem{lewis1985optimal}
F.~L. Lewis, ``Optimal control for singular systems,'' in \emph{IEEE 24th CDC,
  1985}, vol.~24, 1985, pp. 266--272.

\bibitem{li2017network}
C.~Li, H.-D. Chiang, and Z.~Du, ``Network-preserving sensitivity-based
  generation rescheduling for suppressing power system oscillations,''
  \emph{IEEE Trans. on Power Sys.}, vol.~32, no.~5, pp. 3824--3832, 2017.

\bibitem{chen1998linear}
C.-T. Chen, \emph{Linear system theory and design}.\hskip 1em plus 0.5em minus
  0.4em\relax Oxford University Press, Inc., 1998.

\bibitem{vidyasagar2002nonlinear}
M.~Vidyasagar, \emph{Nonlinear Systems Analysis}, 2nd~ed.\hskip 1em plus 0.5em
  minus 0.4em\relax Society for Industrial and Applied Mathematics, 2002.

\bibitem{aminzarey2014contraction}
Z.~Aminzarey and E.~D. Sontagy, ``Contraction methods for nonlinear systems: A
  brief introduction and some open problems,'' in \emph{53rd IEEE CDC}, Dec
  2014, pp. 3835--3847.

\bibitem{smith1992some}
R.~L. Smith, ``Some interlacing properties of the schur complement of a
  hermitian matrix,'' \emph{Linear algebra and its applications}, vol. 177, pp.
  137--144, 1992.

\bibitem{nguyencons18}
H.~D. Nguyen, K.~Dvijotham, and K.~Turitsyn, ``Constructing convex inner
  approximations of steady-state security regions,'' \emph{IEEE Transactions on
  Power Systems}, pp. 1--1, 2018.

\bibitem{nguyenframe18}
H.~D. Nguyen, K.~Dvijotham, S.~Yu, and K.~Turitsyn, ``A framework for robust
  long-term voltage stability of distribution systems,'' \emph{IEEE
  Transactions on Smart Grid}, pp. 1--1, 2018.

\bibitem{li2015tight}
Q.~Li, ``A tight sdp relaxation for miqcqp problems in power systems based on
  disjunctive programming,'' \emph{arXiv:1509.05141}, 2015.

\bibitem{molzahn2018towards}
D.~K. Molzahn and L.~A. Roald, ``Towards an ac optimal power flow algorithm
  with robust feasibility guarantees,'' in \emph{2018 PSCC}.\hskip 1em plus
  0.5em minus 0.4em\relax IEEE, 2018, pp. 1--7.

\bibitem{golestaneh2018ellipsoidal}
F.~Golestaneh, P.~Pinson, R.~Azizipanah-Abarghooee, and H.~B. Gooi,
  ``Ellipsoidal prediction regions for multivariate uncertainty
  characterization,'' \emph{IEEE Transactions on Power Systems}, 2018.

\bibitem{Jwang18}
B.~Chen, C.~Chen, J.~Wang, and K.~L. Butler-Purry, ``Multi-time step service
  restoration for advanced distribution systems and microgrids,'' \emph{IEEE
  Trans. on Smart Grid}, vol.~9, no.~6, pp. 6793--6805, Nov 2018.

\bibitem{toker1995np}
O.~Toker and H.~Ozbay, ``On the np-hardness of solving bilinear matrix
  inequalities and simultaneous stabilization with static output feedback,'' in
  \emph{ACC, Proc. of 1995}, vol.~4.\hskip 1em plus 0.5em minus 0.4em\relax
  IEEE, 1995, pp. 2525--2526.

\bibitem{PSAT}
\BIBentryALTinterwordspacing
F.~Milano, ``{Power System Analysis Toolbox Reference Manual for PSAT},'' 5
  2010. [Online]. Available: \url{http://faraday1.ucd.ie/psat.html.}
\BIBentrySTDinterwordspacing

\bibitem{lofberg2012}
J.~L{\"o}fberg, ``Automatic robust convex programming,'' \emph{Optimization
  methods and software}, vol.~27, no.~1, pp. 115--129, 2012.

\end{thebibliography}

\end{document}